\documentclass[12pt,twoside,reqno]{amsart}

\usepackage{amsmath,amssymb,amsthm, amsfonts }
\usepackage{fancyhdr}
\usepackage{epsfig,color,colordvi}

\usepackage[numbers,sort&compress,square,comma]{natbib}

\theoremstyle{definition}

\theoremstyle{remark}

\theoremstyle{plain}

\newcommand{\be}{\begin{equation}}
\newcommand{\ee}{\end{equation}}

\def\LL{\mathbf{L}}

\def\ee{\mathbb{E}}

\def\rr{\mathbb{R}}
\def\pp{\mathbb{P}}

\def\LL{\mathcal L}

\newtheorem{thm}{Theorem}[section]
\newtheorem{cor}[thm]{Corollary}
\newtheorem{lem}[thm]{Lemma}
\newtheorem{prop}[thm]{Proposition}
\newtheorem{remarks}[thm]{Remarks}

\def\bdes{\begin{description}}
\def\ndes{\end{description}}
\def\beq{\begin{equation}}
\def\deq{\end{equation}}

\def\bdef{\begin{defn}}
\def\ndef{\end{defn}}

\def\bthm{\begin{thm}}
\def\nthm{\end{thm}}

\def\bprop{\begin{prop}}
\def\nprop{\end{prop}}

\def\brmk{\begin{remarks}}
\def\nrmk{\end{remarks}}

\def\bexa{\begin{exa}}
\def\nexa{\end{exa}}

\def\blem{\begin{lem}}
\def\nlem{\end{lem}}

\def\bcor{\begin{cor}}
\def\ncor{\end{cor}}

\def\bexe{\begin{exe}}
\def\nexe{\end{exe}}

\def\bprf{\begin{proof}}
\def\nprf{\end{proof}}

\def\fQ{{{\rm Q}\kern-.65em {}^{{}_/ }\,}}
\def\fQQ{ {{\rm Q}\kern-.57em \scriptscriptstyle{}^{]\kern.055em[}\,}}

\def\ord{\kern0.1em o\kern-0.02em{}_{\ds\breve{}}\kern0.1em}
\def\Ord{\kern0.1em O\kern-0.02em{\ds\breve{}}\kern0.1em}
\def\ds{\displaystyle}

\def\fmonth{\ifcase\month\or Jan\or Feb\or Mar\or Apr
\or May\or Jun\or Jul\or Aug\or Sep
\or Oct\or Nov\or Dec\fi\ }
\def\mmddyyyy{\the\month.\the\day.\the\year}
\def\ddmmyyyy{\the\day.\the\month.\the\year}
\def\Mddyyyy{\fmonth~\the\day,~\the\year}

\providecommand{\pp}[1]{\langle#1\rangle}

\numberwithin{equation}{section}

\allowdisplaybreaks[1]

\begin{document}


\author{Yu Miao}   
\address{College of Mathematics and Information Science, Henan
Normal University, 453007 Henan, China and Department of Mathematics
and Statistics, Wuhan University, 430072 Hubei, China.}
\email{yumiao728@yahoo.com.cn}







\begin{abstract}
 In this paper, we give the central limit theorem and almost sure
central limit theorem for products of some partial sums of
independent identically distributed random variables.

\end{abstract}




\title[CLT and ASCLT for product of some partial sums]
{Central limit theorem and almost sure central limit theorem for the
product of some partial sums}

\maketitle


\def\X{\mathcal X}
\def\B{\mathcal B}

\section{Introduction}
Let $(X_n)_{n\ge 1}$ be a sequence of independent identically
distributed (i.i.d.) positive random variables (r.v.). Recently
there have been several studies to the products of partial sums.
 It is well known that the products of i.i.d. positive, square integrable random variables are
asymptotically log-normal. This fact is an immediate consequence of
the classical central limit theorem (CLT). This point, up to the
knowledge of the author, was first argued by Arnold and Villase\~nr
\cite{AV}, who considered the limiting properties of the sums of
records. In their paper Arnold and Villase\~nr obtained the
following version of the CLT for a sequence of i.i.d. exponential
r.v.'s $(X_n)_{n\ge 1}$ with the mean equal to one:
$$
\frac{\sum_{k=1}^n\log S_k-n\log
n+n}{\sqrt{2n}}\xrightarrow{\LL}\Phi, \ {\rm as}\ n\to\infty,
$$
where $S_k=\sum_{j=1}^kX_j$, $1\le k\le n$, and $\Phi$ is a standard
normal r.v.. Rempa{\l}a and Weso\l owski \cite{RW} have noted that
this limit behavior of a product of partial sums has a universal
character and holds for any sequence of square integrable, positive
i.i.d. random variables. Namely, they have proved the following.

\vskip5pt

\noindent \textbf{Theorem RW} Let $(X_n)_{n\ge 1}$ be a sequence of
i.i.d. positive square integrable random variables with $\ee
X_1=\mu$, $Var X_1=\sigma^2>0$ and the coefficient of variation
$\gamma=\sigma/\mu$. Then
 \beq\label{RM}
\left(\frac{\prod_{k=1}^n S_k}{n!\mu^n}\right)^{1/(\gamma\sqrt{n})}
\xrightarrow{\LL} e^{\sqrt{2}\Phi}.
 \deq
\vskip5pt

Recently, Gonchigdanzan and Rempa{\l}a \cite{GR} discussed an almost
sure limit theorem for the product of the partial sums of i.i.d.
positive random variables as follows.

\vskip5pt

\noindent \textbf{Theorem GR} Let $(X_n)_{n\ge 1}$ be a sequence of
i.i.d. positive square integrable random variables with $\ee
X_1=\mu>0$, $Var X_1=\sigma^2$. Denote $\gamma=\sigma/\mu$ the
coefficient of variation. Then for any real $x$,
 \beq\label{GR}
\lim_{N\to\infty}\frac{1}{\log N}\sum_{n=1}^N\frac 1n
I\left(\left(\frac{\prod_{k=1}^{n}S_k}{n!\mu^n}\right)^{1/(\gamma\sqrt
n)}\le x\right)=F(x), \ \ a.s.
 \deq
where $F(\cdot)$ is the distribution function of the r.v. $e^{\sqrt
2 \Phi}$.

\vskip5pt

 For further discussions of the CLT, the author refers to
\cite{LQ,Qi}. Zhang and Huang \cite{ZH} obtained the invariance
principle of the product of sums of random variables. It is perhaps
worth to notice that by the strong law of large numbers and the
property of the geometric mean it follows directly that
 \beq\label{102}
\left(\frac{\prod_{k=1}^n S_k}{n!}\right)^{1/n} \xrightarrow{a.s.}
\mu
 \deq
if only the existence of the first moment is assumed.

\vskip5pt

Throughout the present paper let $S_{n,k}=\sum_{i=1}^nX_i-X_k$ for
all $n\ge 1, 1\le k\le n$ and we are interested in the similar
results as (\ref{RM}) and (\ref{GR}).

\section{Central limit theorem}

\begin{thm}\label{thm1}
Let $(X_n)_{n\ge 1}$ be a sequence of i.i.d. positive square
integrable random variables with $\ee X_1=\mu$, $Var X_1=\sigma^2>0$
and the coefficient of variation $\gamma=\sigma/\mu$. Then
 \beq\label{thm1-1}
  \left(\frac{\prod_{k=1}^n
S_{n,k}}{(n-1)^n\mu^n}\right)^{1/(\gamma\sqrt{n})} \xrightarrow{\LL}
e^{\Phi}.
 \deq
 where $\Phi$ is a standard
normal r.v.
\end{thm}
\begin{proof}
Let $Y_{i}=(X_i-\mu)/\sigma$, $i=1,2,\cdots$, then
 \beq\label{thm1-00}
 \frac{1}{\gamma\sqrt{n}}\sum_{k=1}^n\left(\frac{S_{n,k}}{(n-1)\mu}-1\right)
 =\frac{1}{\sqrt{n}}\sum_{k=1}^n\left(\frac{\sum_{i\ne k,i\le n}(X_{i}-\mu)}{(n-1)\sigma}\right)
 =\frac{1}{\sqrt{n}}\sum_{k=1}^nY_k.
 \deq
Therefore from the classical central limit theorem and $\ee Y_i=0$,
$Var(Y_i)=1$ for all $i=1, 2,\cdots$, we know that
 \beq\label{thm1-01}
\frac{1}{\gamma\sqrt{n}}\sum_{k=1}^n\left(\frac{S_{n,k}}{(n-1)\mu}-1\right)\xrightarrow{\LL}\Phi.
 \deq
 Furthermore let $C_{n,k}=S_{n,k}/((n-1)\mu)$, $k=1,2,\cdots$. By
 the  strong law of large numbers it follows that for any
 $\delta>0$, $\exists$ $R$ such that,
 $$
 \pp\left(\sup_{n\ge R, 1\le k\le
 n}|C_{n,k}-1|\ge \delta\right)<\delta.
 $$
 Taking $\delta<1/2$, for any $x\in\rr$, we have
$$
\aligned
 \pp\left(\frac{1}{\gamma\sqrt{n}}\sum_{k=1}^n\log(C_{n,k})\le
 x\right)=&\pp\left(\frac{1}{\gamma\sqrt{n}}\sum_{k=1}^n\log(C_{n,k})\le
 x, \sup_{n\ge R,1\le k\le
 n}|C_{n,k}-1|\ge \delta\right)\\
 &\ \ + \pp\left(\frac{1}{\gamma\sqrt{n}}\sum_{k=1}^n\log(C_{n,k})\le
 x, \sup_{n\ge R,1\le k\le
 n}|C_{n,k}-1|< \delta\right)\\
 :=& A_{n}+B_{n}
\endaligned
$$
and
 \beq\label{A}
A_{n}\le \delta.
 \deq
 Next we will control the term $B_{n}$.
By the following logarithm:
$$
 \log(1+x)=x+\frac{x^2}{(1+\theta x)^2}, \
$$
where $\theta\in(0,1)$ depends on $x\in(-1,1)$, we have
$$
\aligned
 &B_{n }=\pp\left(\frac{1}{\gamma\sqrt{n}}\sum_{k=1}^n\log(C_{n,k})\le
 x, \sup_{n\ge R, 1\le k\le
 n}|C_{n,k}-1|< \delta \right)\\
  =&\pp\Big\{\frac{1}{\gamma\sqrt{n}}\sum_{k=1}^{n}(C_{n,k}-1)
  +\frac{1}{\gamma\sqrt{n}}\sum_{k=1}^{n}\frac{(C_{n,k}-1)^2}{(1+\theta_k(C_{n,k}-1))^2}\le x,
 \sup_{n\ge R,1\le k\le
 n }|C_{n,k}-1|< \delta\Big\}\\
 =&\pp\Big\{\frac{1}{\gamma\sqrt{n}}\sum_{k=1}^{n}(C_{n,k}-1)
 +\Big[\frac{1}{\gamma\sqrt{n}}\sum_{k=1}^{n}\frac{(C_{n,k}-1)^2}{(1+\theta_k(C_{n,k}-1))^2}\Big]
 I(\sup_{n\ge R,1\le k\le
 n}|C_{n,k}-1|< \delta)\le
 x\Big\}\\
 & -\pp\Big\{\frac{1}{\gamma\sqrt{n}}\sum_{k=1}^{n}(C_{n,k}-1)\le x,
 \sup_{n\ge R,1\le k\le
 n }|C_{n,k}-1|\ge\delta \Big\}\\
 :=& D_{n }+F_{n },
\endaligned
$$
where $\theta_k, k=1,\cdots,n$ are $(0,1)$-valued random variables
and $F_{n}\le \delta $. To estimate the term $D_{n}$, by the
following elementary inequality: for $|x|<1/2$ and any
$\theta\in(0,1)$ it follows that $x^2/(1+\theta x)^2\le 4x^2$, then
we have
 \beq\label{D}
 \aligned
&\Big[\frac{1}{\gamma\sqrt{n}}\sum_{k=1}^{n}\frac{(C_{n,k}-1)^2}{(1+\theta_k(C_{n,k}-1))^2}\Big]
 I(\sup_{n\ge R,1\le k\le
 n }|C_{n,k}-1|< \delta )
 \le &
 \frac{4}{\gamma\sqrt{n}}\sum_{k=1}^{n}(C_{n,k}-1)^2\xrightarrow{\pp}
 0,
 \endaligned
 \deq
as $n\to\infty$. Relation (\ref{D}) is a consequence of the Markov
inequality, since for any $r>0$,
$$
\pp\left(\frac{4}{\gamma\sqrt{n}}\sum_{k=1}^{n}(C_{n,k}-1)^2\ge
r\right)\le \frac{4}{r\gamma\sqrt{n}}\sum_{k=1}^{n}Var(C_{n,k}-1)
=\frac{4}{r\gamma\sqrt{n}}\frac{n\gamma^2}{n-1}\to 0.
$$
Therefore $D_n\xrightarrow{\LL}\Phi(x)$.  For any $x\in\rr$, we have
$$
\pp\left(\log\left(\frac{\prod_{k=1}^n
S_{n,k}}{(n-1)^n\mu^n}\right)^{1/(\gamma\sqrt{n})}\le
x\right)=\pp\left(\frac{1}{\gamma\sqrt{n}}\sum_{k=1}^n\log
C_{n,k}\le x\right)=A_n+D_n+F_n
$$
which implies our result since the above discussions.

\end{proof}

\section{Almost sure central limit
theorem}

In this section we will consider the almost sure central limit
theorem as (\ref{GR}). Starting with Brosamler \cite{Br} and Schatte
\cite{Sc}, in the past decade several authors investigated the a.s.
central limit theorem and related "logarithmic" limit theorems for
partial sums of independent random variables. The simplest form of
the a.s. central limit theorem (Brosamler, \cite{Br}; Schatte,
\cite{Sc}; Lacey and Philipp, \cite{LP}) states that if $X_1, X_2,
\cdots $ are i.i.d. random variables with mean $0$, variance $1$ and
partial sums $S_n=\sum_{i=1}^nX_i$ then
 \beq\label{AS}
\lim_{N\to\infty}\frac{1}{\log
N}\sum_{k=1}^N\frac{1}{k}I\left(\frac{S_k}{\sqrt k}\le
x\right)=\Phi(x)\ \ a.s., \ \forall x,
 \deq
  where $I$ denotes indicator function. Berkes and Cs\'aki \cite{BC} extended this theory and show that not only the central limit theorem, but every weak
limit theorem for independent random variables, subject to minor
technical conditions, has an analogous almost sure version. However
under our model we only need the simplest version of (\ref{AS}).

\begin{thm}\label{thm2}
Let $(X_n)_{n\ge 1}$ be a sequence of i.i.d. positive square
integrable random variables with $\ee X_1=\mu>0$, $Var
X_1=\sigma^2$. Denote $\gamma=\sigma/\mu$ the coefficient of
variation. Then for any real $x$,
 \beq\label{thm2-1}
\lim_{N\to\infty}\frac{1}{\log N}\sum_{n=1}^N\frac 1n
I\left(\left(\frac{\prod_{k=1}^{n}S_{n,k}}{(n-1)^n\mu^n}\right)^{1/(\gamma\sqrt
n)}\le x\right)=F(x), \ \ a.s.
 \deq
where $F(\cdot)$ is the distribution function of the r.v.
$e^{\Phi}$.
\end{thm}
\begin{proof}
Let $Y_{i}=(X_i-\mu)/\sigma$, $i=1,2,\cdots$, then $\ee Y_i=0$ and
$Var(Y_i)=1$ for all $i=1,2,\cdots$ and from (\ref{thm1-00}), for
any real $x$, we have
 \beq\label{thm2-2}
\aligned
 &\lim_{N\to\infty}\frac{1}{\log N}\sum_{n=1}^N\frac 1n
I\left(
\frac{1}{\gamma\sqrt{n}}\sum_{k=1}^n\left(\frac{S_{n,k}}{(n-1)\mu}-1\right)\le
x\right)\\
  =&\lim_{N\to\infty}\frac{1}{\log N}\sum_{n=1}^N\frac 1n
I\left(\frac{1}{\sqrt{n}}\sum_{k=1}^nY_k\le x\right)=\Phi(x)\ \ \
a.s.
\endaligned
 \deq
 Note that in order to prove (\ref{thm2-1}) it is sufficient to
show that for any real $x$,
 \beq\label{thm2-3}
\lim_{N\to\infty}\frac{1}{\log N}\sum_{n=1}^N\frac 1n
I\left(\frac{1}{\gamma\sqrt
n}\sum_{k=1}^{n}\log\frac{S_{n,k}}{(n-1)\mu}\le x\right)=\Phi(x), \
\ a.s.
 \deq
To this end let, as before, $C_{n,k}=S_{n,k}/((n-1)\mu)$ and note
that by the law of the iterated logarithm we have for $n\to\infty$,
$$
\max_{1\le k\le n}|C_{n,k}-1|=O\left(\left(\frac{\log\log
n}{n}\right)^{1/2}\right)\ \ \ a.s.
$$
Since for $|x|<1$ we have $\log (1+x)=x+R(x)$ with $\lim_{x\to
0}R(x)/x^2=1/2$, thus
$$
\aligned
 & \left|\sum_{k=1}^n\log
C_{n,k}-\sum_{k=1}^n(C_{n,k}-1)\right|\ll\sum_{k=1}^n(C_{n,k}-1)^2\\
&\ll\sum_{k=1}^n\frac{\log\log n}{n}\ll\log\log n\log n\ \ a.s.
\endaligned
$$
where $"\ll"$ denote the inequality $"\le"$ up to some universal
constant. Hence for almost every $\omega$ and any $\varepsilon>0$
there exists $n_0=n_0(\omega,\varepsilon,x)$ such that for $n\ge
n_0$
$$
\aligned I\left(
\frac{1}{\gamma\sqrt{n}}\sum_{k=1}^n\left(\frac{S_{n,k}}{(n-1)\mu}-1\right)\le
x-\varepsilon\right)\le& I\left(\frac{1}{\gamma\sqrt
n}\sum_{k=1}^{n}\log\frac{S_{n,k}}{(n-1)\mu}\le x\right)\\
\le & I\left(
\frac{1}{\gamma\sqrt{n}}\sum_{k=1}^n\left(\frac{S_{n,k}}{(n-1)\mu}-1\right)\le
x+\varepsilon\right)
\endaligned
$$
thus (\ref{thm2-2}) implies (\ref{thm2-3}).
\end{proof}

\section*{Acknowledgements}
Here the author is very grateful to the referee for his valuable
report and some helpful suggestions.

\end{document}